\documentclass[12pt,twoside]{article}

\usepackage{amsfonts}
\usepackage{amssymb}
\usepackage{amsmath}
\usepackage{graphicx}
\usepackage{color}

\newcommand{\ignore}[1]{}

\def\@begintheorem#1#2{\par\bgroup{\sc #1\ #2. }\it\ignorespaces}
\def\@opargbegintheorem#1#2#3{\par\bgroup{\sc #1\ #2\ (#3). }\it\ignorespaces}
\def\@endtheorem{\egroup}
\newtheorem{theorem}{Theorem}[section]
\newtheorem{corollary}[theorem]{Corollary}
\newtheorem{lemma}[theorem]{Lemma}
\newtheorem{proposition}[theorem]{Proposition}
\newtheorem{example}[theorem]{Example}
\newtheorem{algorithm}[theorem]{Algorithm}
\newtheorem{definition}[theorem]{Definition}
\newtheorem{remark}[theorem]{Remark}
\newcommand{\bt}[1]{\begin{theorem}\label{#1}}
\newcommand{\bc}[1]{\begin{corollary}\label{#1}}
\newcommand{\bl}[1]{\begin{lemma}\label{#1}}
\newcommand{\bp}[1]{\begin{proposition}\label{#1}}
\newcommand{\be}[1]{\begin{example}\rm\label{#1}}
\newcommand{\ba}[1]{\begin{algorithm}\rm\label{#1}}
\newcommand{\bd}[1]{\begin{definition}\rm\label{#1}}
\newcommand{\br}[1]{\begin{remark}\rm\label{#1}}

\newcommand{\et}{\end{theorem}}
\newcommand{\ec}{\end{corollary}}
\newcommand{\el}{\end{lemma}}
\newcommand{\ep}{\end{proposition}}
\newcommand{\ee}{\end{example}}
\newcommand{\ea}{\end{algorithm}}
\newcommand{\ed}{\end{definition}}
\newcommand{\er}{\end{remark}}

\def\R{\mathbb{R}}
\def\Z{\mathbb{Z}}

\def \G {{\cal G}}
\def \l {\langle}
\def \r {\rangle}

\def \G {{\cal G}}

\def \l {\langle}
\def \r {\rangle}

\def \A {A^{(n)}}

\usepackage{enumerate}
\def\ve#1{\mathchoice{\mbox{\boldmath$\displaystyle\bf#1$}}
{\mbox{\boldmath$\textstyle\bf#1$}}
{\mbox{\boldmath$\scriptstyle\bf#1$}}
{\mbox{\boldmath$\scriptscriptstyle\bf#1$}}}

\newcommand\veb{{\ve b}}
\newcommand\vecc{{\ve c}}

\newcommand\veg{{\ve g}}
\newcommand\veh{{\ve h}}
\newcommand\vel{{\ve l}}

\newcommand\veu{{\ve u}}
\newcommand\vev{{\ve v}}
\newcommand\vew{{\ve w}}
\newcommand\vex{{\ve x}}
\newcommand\vey{{\ve y}}
\newcommand\vez{{\ve z}}

\def \Graver {{\G}}

\newcommand{\boproof}{\textbf{Proof.} }
\newcommand{\eoproof}{\hspace*{\fill} $\square$ \vspace{5pt}}

\begin{document}

\title{\bf N-fold integer programming in cubic time}

\author
{Raymond Hemmecke\thanks{Technische Universit\"at Munich, Germany.}
\and
Shmuel Onn
\thanks{Technion -- Israel Institute of Technology.
Supported in part by the Israel Science Foundation.}
\and
Lyubov Romanchuk
\thanks{Technion -- Israel Institute of Technology.
Supported in part by the Technion Graduate School.}
}

\date{}

\maketitle

\begin{abstract}
N-fold integer programming is a fundamental problem with a variety
of natural applications in operations research and statistics.
Moreover, it is universal and provides a new, variable-dimension,
parametrization of all of integer programming. The fastest algorithm
for $n$-fold integer programming predating the present article runs
in time $O\left(n^{g(A)}L\right)$ with $L$ the binary length of the
numerical part of the input and $g(A)$ the so-called Graver complexity
of the bimatrix $A$ defining the system. In this article we provide a
drastic improvement and establish an algorithm which runs in time
$O\left(n^3 L\right)$ having cubic dependency on $n$ regardless of
the bimatrix $A$. Our algorithm can be extended to separable convex piecewise
affine objectives as well, and also to systems defined by bimatrices
with variable entries.  Moreover, it can be used to define a
hierarchy of approximations for any integer programming problem.
\end{abstract}

\section{Introduction}
\label{i}

N-fold integer programming is the following problem in
variable dimension $nt$,
\begin{equation}\label{NIP}
\min\left\{\vew\vex\ :\ \A\vex=\veb\,,\
\vel\leq\vex\leq\veu\,,\ \vex\in\Z^{nt}\right\}\ ,
\end{equation}

where
\begin{equation}\label{NFold}
A^{(n)}\quad:=\quad
\left(
\begin{array}{cccc}
  A_1    & A_1    & \cdots & A_1    \\
  A_2    & 0      & \cdots & 0      \\
  0      & A_2    & \cdots & 0      \\
  \vdots & \vdots & \ddots & \vdots \\
  0      & 0      & \cdots & A_2    \\
\end{array}
\right)\quad
\end{equation}
is an $(r+ns)\times nt$ matrix which is the {\em $n$-fold product}
of a fixed {\em $(r,s)\times t$ bimatrix} $A$, that is, of a matrix $A$
consisting of two blocks $A_1$, $A_2$, with $A_1$ its $r\times t$
submatrix consisting of the first $r$ rows and $A_2$ its $s\times t$
submatrix consisting of the last $s$ rows.

This is a fundamental problem with a variety of natural applications
in operations research and statistics which, along with extensions
and variations, include multiindex and multicommodity transportation
problems, privacy and disclosure control in statistical databases,
and stochastic integer programming. We briefly discuss some of these
applications in Section \ref{sc} (Corollaries \ref{Multicommodity}
and \ref{Privacy}). For more information see e.g.
\cite{DHORW, DHOW, DFRSZ, HKW, HOW1, HOW2, HSc, KMW, LS, SZP},
\cite[Chapters 4,5]{Onn}, and the references therein.

Moreover, $n$-fold integer programming is universal \cite{DO} and provides
a new, variable-dimension, parametrization of all of integer programming:
every program is an $n$-fold program for some $m$ over the bimatrix $A:=A(m)$
with first block $A_1$ the $3m\times 3m$ identity matrix and second
block $A_2$ the $(3+m)\times 3m$ incidence matrix of the complete bipartite
graph $K_{3,m}$. We make further discussion of this in Section \ref{paah}.

The fastest algorithm for $n$-fold integer programming predating
the present article is in \cite{HOW1} and runs in time
$O\left(n^{g(A)}L\right)$ with $L=\l \vew,\veb,\vel,\veu \r$
the binary length of the numerical part of the input, and $g(A)$
the so-called {\em Graver complexity} of the bimatrix $A$.
Unfortunately, the Graver complexity is typically very large \cite{BO,SS}:
for instance, the bimatrices $A(m)$ mentioned above have Graver
complexity $g(A(m))=\Omega(2^m)$, yielding polynomial but very large
$n^{\Omega(2^m)}$ dependency of the running time on $n$.

In this article we provide a drastic improvement and establish an
algorithm which runs in time $O\left(n^3 L\right)$ having the
cubic dependency on $n$ which is alluded to in the title, regardless
of the fixed bimatrix $A$. So the Graver complexity $g(A)$ now drops
down from the exponent of $n$ to the constant multiplying $n^3$.
This is established in Section \ref{ta} (Theorem \ref{Main}).
Moreover, our construction can be used to define
a natural hierarchy of approximations for (\ref{NIP}) for the
bimatrices $A(m)$ with $m$ variable, and therefore, by the
universality theorem of \cite{DO}, for any integer programming
problem. These approximations are currently under study,
implementation and testing, and will be discussed briefly
in Section \ref{paah} and in more detail elsewhere.

Our algorithm extends, moreover, for certain nonlinear objective
functions: using results of \cite{MSW} on certain optimality criteria,
we provide in Section \ref{etno} an optimality certification procedure
for separable convex objectives whose time complexity is linear
in $n$ (Theorem \ref{Separable}) and an algorithm for solving
problems with separable convex piecewise affine objectives whose
time complexity is again cubic in $n$ (Theorem \ref{Piecewise}).
Furthermore, the algorithm also leads to the first polynomial
time solution of $n$-fold integer programming problems over bimatrices
with variable entries (Theorem \ref{Variable}).

\section{Notation and preliminaries}
\label{nap}

We start with some notation and review of some preliminaries on
Graver bases and $n$-fold integer programming that we need later on.
See the book \cite{Onn} for more details.

Graver bases were introduced in \cite{Gra} as optimality certificates
for integer programming. Define a partial order $\sqsubseteq$ on $\R^n$
by $\vex\sqsubseteq\vey$ if $x_iy_i\geq 0$ and $|x_i|\leq |y_i|$
for all $i$. So $\sqsubseteq$ extends the coordinate-wise partial
order $\leq$ on the nonnegative orthant $\R_+^n$ to all of $\R^n$.
By a classical lemma of Gordan, every subset $Z\subseteq \Z^n$
has finitely-many $\sqsubseteq$-minimal elements, that is, $\vex\in Z$
such that no other $\vey\in Z$ satisfies $\vey\sqsubseteq\vex$.

We have the following fundamental definition from \cite{Gra}.

\bd{GraverBasisDefinition}
The {\em Graver basis} of an integer $m\times n$ matrix $A$
is defined to be the finite set $\G(A)\subset\Z^n$ of
$\sqsubseteq$-minimal elements in $\{\vex\in\Z^n\,:\, A\vex=0,\ \vex\neq\ve 0\}$.
\ed

For instance, the Graver basis of the matrix $A:=(1\,\ 2\,\ 1)$ consists of $8$ vectors,
$$\G(A)\ =\
\pm\left\{\, (2\,\ -1\,\ 0)\,,\ (0\,\ -1\,\ 2)\,,\
(1\,\ 0\,\ -1)\,,\ (1\,\ -1\,\ 1)\, \right\}\ .$$

Consider the general integer programming problem in standard form,
\begin{equation}\label{IP}
\min\left\{\vew\vex\ :\ A\vex=\veb\,,\
\vel\leq\vex\leq\veu\,,\ \vex\in\Z^n\right\}\ .
\end{equation}
A {\em feasible step} for feasible point $\vex$ in (\ref{IP})
is any vector $\vev$ such that $\vex+\vev$ is also feasible,
that is, $A\vev=\ve0$ so $A(\vex+\vev)=\veb$,
and $\vel\leq\vex+\vev\leq\veu$. An {\em augmenting step} for $\vex$
is a feasible step $\vev$ such that $\vex+\vev$ is better,
that is, $\vew\vev<0$ so $\vew(\vex+\vev)<\vew\vex$.

Graver has shown that a feasible point $\vex$ in (\ref{IP})
is optimal if and only if there is no element $\veg\in\G(A)$
in the Graver basis of $A$ which is an augmenting step for $\vex$.

This suggests the following simple augmentation scheme: start from any
feasible point in (\ref{IP}) and iteratively augment it to an optimal
solution using Graver augmenting steps $\veg\in\G(A)$ as long as possible.
While the number of iterations in this simple scheme as is may be
exponential, it was recently shown in \cite{HOW1} that if in each iteration
the best possible augmenting step of the form $\gamma\veg$ with
$\gamma$ positive integer and $\veg\in\Graver(A)$ is taken,
then the number of iterations does become polynomial.
In what follows, we call an augmenting step which is at least as good
as the best possible augmenting step $\gamma\veg$ with
$\veg\in\Graver(A)$, a {\em Graver-best augmenting step}.

It was shown in \cite{DHOW} that for fixed bimatrix $A$, the
Graver basis $\G\left(\A\right)$ of the $n$-fold product of $A$
can be computed in time polynomial in $n$. Thus, to find a Graver-best
augmenting step of the form $\gamma\veg$ for an $n$-fold integer
program (\ref{NIP}), it is possible, as shown in \cite{HOW1},
to check each element $g\in\G\left(\A\right)$, and for each,
find the best possible step size $\gamma$. However, as we
explain below, the Graver basis $\G\left(\A\right)$ is very large.
Therefore, in this article, we do it the other way around. For each of
$O(n)$ critical positive integer potential step sizes $\gamma$, we
determine an augmenting step $\gamma\veh$ which is at least as
good as the best possible augmenting step $\gamma\veg$
with $g\in\G\left(\A\right)$.
We then show that the best among these steps over all such
$\gamma$ is a Graver-best step.

In preparation for this, we need to review some material on Graver bases
of $n$-fold products. Let $A$ be a fixed integer $(r,s)\times t$ bimatrix.
For any $n$ we write each vector $\vex\in\Z^{nt}$ as a tuple
$\vex=(\vex^1,\dots,\vex^n)$ of $n$ {\em bricks} $\vex^i\in\Z^t$.
It has been shown in \cite{AT}, \cite{SS}, and \cite{HS}, in increasing
generality, that for every bimatrix $A$, the number of nonzero bricks
appearing in any element in the Graver basis $\G\left(\A\right)$
for any $n$ is bounded by a constant independent of $n$.
So we can make the following definition.

\bd{GraverComplexity}
The {\em Graver complexity} of an integer bimatrix $A$ is
defined to be the largest number $g(A)$ of nonzero bricks
$\veg^i$ in any element $\veg\in\G\left(\A\right)$ for any $n$.
\ed

This was used in \cite{DHOW} to show that the Graver basis
$\G(A^{(n)})$ has a polynomial number $O(n^{g(A)})$ of elements
and is computable in time $O(n^{g(A)})$ polynomial in $n$.
Thus, the computation of a Graver-best augmenting step in
\cite{HOW1} was done by finding the best step size
$\gamma$ for each of these $O(n^{g(A)})$ elements of $\G(A^{(n)})$,
resulting in polynomial but very large $O(n^{g(A)})$ dependency
of the running time on $n$. In Section \ref{ta} we
show how to find a Graver-best step without constructing
$\G(A^{(n)})$ explicitly in quadratic time $O(n^2)$
regardless of the bimatrix $A$ and its Graver complexity.

We conclude this subsection with some remarks about
complexity and finiteness. The {\em binary length}
of an integer number $z$ is the number of bits
in its binary encoding, which is $O(\log|z|)$,
and is denoted by $\l z\r$. The binary length $\l\vez\r$
of an integer vector $\vez$ is the sum of binary lengths
of its entries. We denote by $L$ the binary length of
all numerical part of the input. In particular
$L=\l \vew,\veb,\vel,\veu \r$ for problem (\ref{NIP}).
All numbers manipulated by our algorithms remain polynomial
in the binary length of the input and our algorithms are
polynomial time in the Turing machine model. But we are
mostly interested in the number of arithmetic operations
performed (additions, multiplications, divisions, comparisons),
so time in our complexity statements is the number of such
operations as in the real arithmetic model of computation.

For simplicity of presentation we assume throughout that
all entries of the bounds $\vel,\veu$ in program (\ref{NIP})
are finite and hence the set of feasible points in (\ref{NIP})
is finite. This is no loss of generality since, as is well known,
it is always possible to add suitable polynomial upper and
lower bounds without excluding some optimal solution if any.

\section{The algorithm}
\label{ta}

We now show how to decide if a given feasible point
$\vex=(\vex^1,\dots,\vex^n)$ in (\ref{NIP}) is optimal
in linear time $O(n)$, and if not, determine a Graver-best step
$\gamma\veg$ for $\vex$ in quadratic time $O(n^2)$. This is then
incorporated into an iterative algorithm for solving (\ref{NIP}).

We begin with a lemma about elements of Graver bases of $n$-fold products.

\bl{StateSpace}
Let $A$ be integer $(r,s)\times t$ bimatrix with Graver complexity $g(A)$. Let
\begin{eqnarray}\label{Space}
Z(A)\ :=\ \left\{\vez\in\Z^t\,:\, \vez\
\mbox{is the sum of at most $g(A)$ elements of $\G(A_2)$}\right\}\ .
\end{eqnarray}
Then for any $n$, any $\veg\in\G\left(\A\right)$
and any $I\subseteq\{1,\dots,n\}$, we have $\sum_{i\in I}\veg^i\in Z(A)$.
\el
\boproof
Consider any Graver basis element $\veg\in\G\left(\A\right)$ for some $n$.
Then $\A\veg=\ve0$ and hence $\sum_{i=1}^n A_1\veg^i=\ve0$ and
$A_2\veg^i=\ve0$ for all $i$. Therefore (see \cite[Chapter 4]{Onn})
each $\veg^i$ can be written as the sum $\veg^i=\sum_{j=1}^{k_i}\veh^{i,j}$
of some elements $\veh^{i,j}\in\G(A_2)$ for all $i,j$.
Let $m:=k_1+\cdots+k_n$ and let $\veh$ be the vector
$$\veh\ :=\ (\veh^{1,1},\dots,\veh^{1,k_1},\dots,
\veh^{n,1},\dots,\veh^{n,k_n})\ \in\ \Z^{mt}\ .$$
Then $\sum_{i,j}A_1\veh^{i,j}=\ve0$ and $A_2\veh^{i,j}=\ve0$
for all $i,j$ and hence $A^{(m)}\veh=\ve0$. We claim that
moreover, $\veh\in\G\left(A^{(m)}\right)$. Suppose indirectly
this is not the case. Then there is an $\bar\veh\in\G\left(A^{(m)}\right)$
with $\bar\veh\sqsubset\veh$. But then the vector $\bar\veg\in\Z^{nt}$
defined by $\bar\veg^i:=\sum_{j=1}^{k_i}\bar\veh^{i,j}$ for all $i$
satisfies $\bar\veg\sqsubset\bar\veg$ contradicting $\veg\in\G\left(\A\right)$.
This proves the claim. Therefore, by Definition \ref{GraverComplexity}
of Graver complexity, the number of nonzero bricks $\veh^{i,j}$ of
$\veh$ is at most $g(A)$. So for every $I\subseteq\{1,\dots,n\}$,
we have that $\sum_{i\in I}\veg^i=\sum_{i\in I}\sum_{j=1}^{k_i}\veh^{i,j}$
is a sum of at most $g(A)$ nonzero elements $\veh^{i,j}\in\G(A_2)$
and hence $\sum_{i\in I}\veg^i\in Z(A)$.
\eoproof

\be{Example}
Let $A:=A(3)$ be the $(9,6)\times 9$ bimatrix mentioned in the  introduction,
which arises in the universality of $n$-fold integer programming
discussed further in Section \ref{paah}, having first block $A_1=I_9$
the $9\times 9$ identity matrix and second block the following
$6\times 9$ incidence matrix of the complete bipartite graph $K_{3,3}$,
\begin{equation*}\label{matrix}
A_2\ =\
\left(
\begin{array}{ccccccccc}
  1 & 0 & 0 & 1 & 0 & 0 & 1 & 0 & 0 \\
  0 & 1 & 0 & 0 & 1 & 0 & 0 & 1 & 0 \\
  0 & 0 & 1 & 0 & 0 & 1 & 0 & 0 & 1 \\
  1 & 1 & 1 & 0 & 0 & 0 & 0 & 0 & 0 \\
  0 & 0 & 0 & 1 & 1 & 1 & 0 & 0 & 0 \\
  0 & 0 & 0 & 0 & 0 & 0 & 1 & 1 & 1 \\
\end{array}
\right)\ .
\end{equation*}
Since $A_2$ is totally unimodular, its Graver basis $\G(A_2)$
consists of the $30$ vectors in $\{0,\pm1\}^9$ supported
on circuits of $K_{3,3}$ with alternating $\pm1$, see \cite{Onn}.
Also, it is known that the Graver complexity of this bimatrix is $g(A)=9$,
see \cite{BO,SS}. Therefore, the set $Z(A)$ in (\ref{Space}) which corresponds
to $A$, consists of all sums of at most $9$ such circuit
vectors, and turns out to be comprised of $42931$ vectors in $\Z^9$, such as
$$\left(\begin{array}{ccccccccc}
9 & -2 & -7 & -4 & 5 & -1 & -5 & -3 & 8 \end{array}\right)\ .$$
\ee

We now define a dynamic program, that is, a weighted digraph, which
will enable to find a Graver-best step $\gamma\veg$ for a feasible
point $\vex$ of (\ref{NIP}) or detect that none exists.

\bd{DynamicProgram}{\bf (the dynamic program)}
Let $A$ be a fixed $(r,s)\times t$ bimatrix and let $g(A)$ be its Graver complexity.
Given $n$, $\vew,\veb,\vel,\veu$, feasible point $\vex$ in (\ref{NIP}),
and positive integer $\gamma$, define a weighted digraph as follows.
Its vertices are partitioned into $n+1$ stages defined
in terms of the fixed finite set $Z(A)\subset\Z^t$ in (\ref{Space}), by
$$
S_0:=\{\ve0\}\,,\quad S_1:=S_2:=\ \cdots\ :=S_{n-1}:=Z(A)
\,, \quad S_n:=\{\vez\in Z(A)\,:\,A_1\vez=\ve0\}\ .
$$
Denote the vertices of $S_i$ by $\veh^i\in\Z^t$. Introduce an arc $(\veh^{i-1},\veh^i)$
from $\veh^{i-1}\in S_{i-1}$ to $\veh^i\in S_i$ if $\veg^i:=\veh^i-\veh^{i-1}\in Z(A)$
and $\vel^i\leq\vex^i+\gamma\veg^i\leq\veu^i$, and give it weight $\vew^i\veg^i$.
\ed

To each dipath $\veh=(\veh^0,\veh^1,\dots,\veh^n)$ from $S_0$ to $S_n$ in this digraph
we associate a vector $\veg(\veh):=(\veh^1-\veh^0,\dots,\veh^n-\veh^{n-1})\in\Z^{nt}$.
Note that $\ve0\in Z(A)$ and hence the trivial path $\veh=(\ve0,\dots,\ve0)$
with weight $0$ and vector $\veg(\veh)=(\ve0,\dots,\ve0)$ always exists.
Note also that $\vew\veg(\veh)=\sum_{i=1}^n\vew^i(\veh^i-\veh^{i-1})$
is precisely the weight of the dipath $\veh$.

The following lemma relates this dynamic program to Graver augmentations.
\bl{AugmentingStep}
A feasible step $\gamma\veg$ for $\vex$ which satisfies
$\vew(\vex+\gamma\veg)\leq\vew(\vex+\gamma\bar\veg)$
for any feasible step $\gamma\bar\veg$ with $\bar\veg\in\G\left(\A\right)$
can be constructed in linear time $O(n)$.
\el

\boproof
Let $\veh$ be a minimum weight dipath from $S_0$ to $S_n$ and
let $\veg:=\veg(\veh)$ be the vector associated with $\veh$.
We claim that $\gamma\veg$ is the desired feasible step for $\vex$.

We begin with the complexity statement.
Since $A$ is fixed, so is $g(A)$, and hence so is each $S_i$.
As the digraph is acyclic, the minimum weight dipath from $S_0$
to $\veh^i\in S_i$ decomposes into a minimum weight dipath from $S_0$ to some
$\veh^{i-1}\in S_{i-1}$ plus the arc from $\veh^{i-1}$ to $\veh^i$. Thus, we have
to check at most a constant number $|S_{i-1}|\cdot|S_i|$ of such pairs
$(\veh^{i-1},\veh^i)$ to find the minimum weight dipaths from $S_0$ to every
$\veh^i\in S_i$ given the minimum weight dipaths from $S_0$ to every
$\veh^{i-1}\in S_{i-1}$. Repeating this for each of the sets
$S_1,\dots,S_n$ one after the other takes $O(n)$ time.

We next show that $\gamma\veg$ is a good feasible step.
Since $\veg^i=\veh^i-\veh^{i-1}$ and $(\veh^{i-1},\veh^i)$ is an arc, we have
$\vel^i\leq\vex^i+\gamma\veg^i\leq\veu^i$, and $\veg^i\in Z(A)$
and hence $A_2\veg^i=\ve0$, for all $i$. Also, $\sum_{i=1}^n\veg^i=\veh^n\in S_n$
and hence $\sum_{i=1}^n A_1\veg^i= A_1\veh^n=\ve0$. So $\vex+\gamma\veg$
is feasible in (\ref{NIP}). Moreover, $\vew\veg=\sum_{i=1}^n\vew^i\veg^i$
is the weight of the minimum weight dipath $\veh$.
Now consider any feasible step $\gamma\bar\veg$ for $\vex$ with
$\bar\veg\in\Graver\left(\A\right)$. Define $\bar\veh^i:=\sum_{j\leq i}\bar\veg^j$
for all $i$. Then $\bar\veh^i\in Z(A)$ for all $i$ by Lemma \ref{StateSpace}.
Moreover, $A_1\bar\veh^n=A_1\sum_{j=1}^n\bar\veg^j=\ve0$.
Therefore $\bar\veh^i\in S_i$ for all $i$.
Furthermore, $\bar\veh^i-\bar\veh^{i-1}=\bar\veg^i\in Z(A)$
and $\vel^i\leq\vex^i+\gamma\bar\veg^i\leq\veu^i$ and therefore
$(\bar\veh^{i-1},\bar\veh^i)$ is an arc of weight $\vew^i\bar\veg^i$
for all $i$. So $\veh=(\veh^0,\veh^1,\dots,\veh^n)$ is a dipath
from $S_0$ to $S_n$ with weight $\vew\bar\veg=\sum_{i=1}^n\vew^i\bar\veg^i$
and associated vector $\veg(\bar\veh)=\bar\veg$. Since $\veh$ is a minimum
weight dipath, $\vew\veg\leq\vew\bar\veg$ and so
$\vew(\vex+\gamma\veg)\leq \vew(\vex+\gamma\bar\veg)$.
\eoproof

\br{LinearCertification}{\bf (optimality certification in linear time)}
As noted in Section \ref{nap}, a feasible point $\vex$ in integer
program (\ref{NIP}) is optimal if and only if there is no Graver
augmenting step $\veg\in\G\left(\A\right)$ for $\vex$. Thus, with
$\gamma:=1$, Lemma \ref{AugmentingStep} implies that the optimality of a
feasible point $\vex$ in (\ref{NIP}) can be determined in linear time $O(n)$.
\er

The next lemma shows that we can quickly find a Graver-best augmentation.

\bl{BestGraverStep}
A feasible step $\gamma\veg$ for $\vex$ satisfying
$\vew(\vex+\gamma\veg)\leq \vew(\vex+\bar\gamma\bar\veg)$
for any feasible step $\bar\gamma\bar\veg$ with $\bar\gamma\in\Z_+$ and
$\bar\veg\in\G\left(\A\right)$ can be found in quadratic time $O(n^2)$.
\el

\boproof
If $Z(A)=\{\ve0\}$ then $\G\left(\A\right)=\emptyset$ by
Lemma \ref{StateSpace}, so $\gamma \veg:=\ve0$ will do.
Otherwise, construct a set $\Gamma$ of $O(n)$ positive integers
in $O(n)$ time as follows: for every $i=1,\dots,n$ and every
$\vez\in Z(A)\setminus\{\ve0\}$ determine the largest positive integer
$\gamma$ such that $\vel^i\leq\vex^i+\gamma\vez\leq\veu^i$ and include
it in $\Gamma$. Now, for each $\gamma\in\Gamma$, construct and solve the
corresponding dynamic program, resulting in total of $O(n^2)$ time
by Lemma \ref{AugmentingStep}.
Let $\gamma\veg$ be that feasible step for $\vex$ which attains minimum
value $\vew\gamma\veg$ among the best steps obtained from all these
dynamic programs, and let $\bar\gamma\bar\veg$ be that feasible step
for $\vex$ which attains minimum value $\vew\bar\gamma\bar\veg$
among $\bar\gamma\bar\veg$ with $\bar\veg\in\G\left(\A\right)$ if any.
Assume that $\vew\bar\gamma\bar\veg<0$ as otherwise we are done
since $\vew\gamma\veg\leq\vew\gamma\ve0=0$.
Then $\bar\gamma$ is the largest positive integer such that
$\vel\leq\vex+\bar\gamma\bar\veg\leq\veu$ since otherwise
the step $(\bar\gamma+1)\bar\veg$ will be feasible and better. So for
some $i=1,\dots,n$, it must be that $\bar\gamma$ is the largest positive
integer such that $\vel^i\leq\vex^i+\bar\gamma\bar\veg^i\leq\veu^i$.
Since $\bar\veg\in\G\left(\A\right)$, it follows from
Lemma \ref{StateSpace} that $\bar\veg^i\in Z(A)$.
Therefore $\bar\gamma\in\Gamma$. Now let $\bar\gamma\hat\veg$ be
the best step attained from the dynamic program of $\bar\gamma$. Then
$\vew\gamma\veg\leq\vew\bar\gamma\hat\veg$ by choice of $\gamma\veg$ and
$\vew\bar\gamma\hat\veg\leq\vew\bar\gamma\bar\veg$ by Lemma \ref{AugmentingStep}.
Therefore $\vew(\vex+\gamma\veg)\leq \vew(\vex+\bar\gamma\bar\veg)$ as claimed.
\eoproof

We next show, following \cite{HOW1}, that repeatedly applying
Graver-best augmenting steps, we can augment an initial feasible
point for (\ref{NIP}) to an optimal one efficiently.

\bl{Iterations}
For any fixed bimatrix $A$ there is an algorithm that, given
$n$, $\vew,\veb,\vel,\veu$, and feasible point $\vex$ for (\ref{NIP}),
finds an optimal solution $\vex^*$ for (\ref{NIP}) in time $O(n^3 L)$.
\el
\boproof
Iterate the following: find by the algorithm of
Lemma \ref{BestGraverStep} a Graver-best step
$\gamma\veg$ for $\vex$; if it is augmenting then set
$\vex:=\vex+\gamma\veg$ and repeat, else $\vex^*:=\vex$ is optimal.

To bound the number of iterations, following \cite{HOW1},
note that while $\vex$ is not optimal,
and $\vex^*$ is some optimal solution, we have that
$\vex^*-\vex=\sum_{i=1}^k\gamma_i\veg_i$ is a nonnegative
integer combination of Graver basis elements
$\veg_i\in\G\left(\A\right)$ all lying in the same orthant,
and hence each $\vex+\gamma_i\veg_i$ is feasible in (\ref{NIP}).
Moreover, by the integer Carath\'eodory theorem of \cite{CFS,Seb},
we can assume that $k\leq 2(nt-1)$. Letting $\gamma_i\veg_i$ be a
summand attaining minimum $\vew\gamma_i\veg_i$, and letting
$\gamma\veg$ be a Graver-best augmenting step for $\vex$ obtained
from the algorithm of Lemma \ref{BestGraverStep}, we find that
$$\vew(\vex+\gamma\veg)-\vew\vex\ \ \leq\ \
\vew(\vex+\gamma_i\veg_i)-\vew\vex\ \ \leq\ \
{1\over 2(nt-1)}\left(\vew\vex^*-\vew\vex\right)\ .$$
So the Graver-best step provides an improvement which is
a constant fraction of the best possible improvement,
and this can be shown to lead to a bound of $O(nL)$ on the number
of iterations to optimality, see \cite{HOW1} for more details.
Since each iteration takes $O(n^2)$ time by Lemma \ref{BestGraverStep},
the overall running time is $O(n^3 L)$ as claimed.
\eoproof

We next show how to find an initial feasible point for
(\ref{NIP}) with the same complexity. We follow the approach
of \cite{DHOW} using a suitable auxiliary $n$-fold program.

\bl{Feasibility}
For any fixed bimatrix $A$ there is an algorithm that, given
$n$, $\veb,\vel,\veu$, either finds a feasible point $\vex$
for (\ref{NIP}) or asserts that none exists, in time $O(n^3 L)$.
\el
\boproof
Construct an auxiliary $n$-fold integer program
\begin{equation}\label{Auxiliary}
\min\left\{\bar\vew\vez\ :\ {\bar A}^{(n)}\vez=\veb\,,\
\bar\vel\leq\vez\leq\bar\veu\,,\ \vez\in\Z^{n(t+2r+2s)}\right\}
\end{equation}
as follows.
First, construct a new fixed $(r,s)\times(t+2r+2s)$ bimatrix $\bar A$ with
$${\bar A}_1\ :=\
\left(
\begin{array}{ccccc}
  A_1 & I_r & -I_r  & 0_{r\times s} & 0_{r\times s}
\end{array}
\right)\ ,\quad\quad
{\bar A}_2\ :=\
\left(
\begin{array}{ccccc}
  A_2  & 0_{s\times r} & 0_{s\times r}    & I_s & -I_s
\end{array}
\right)\ .
$$
Now, the $n(t+2r+2s)$ variables $\vez$ have a natural partition
into $nt$ original variables $\vex$ and $n(2r+2s)$ new auxiliary
variables $\vey$. Keep the original lower and upper bounds on the
original variables and introduce lower bound $0$ and
upper bound $\|\veb\|_\infty$ on each auxiliary variable.
Let the new objective $\bar\vew\vez$ be the sum of auxiliary
variables. Note that the binary length of the auxiliary program
satisfies $\bar L=O(L)$ and an initial feasible point $\bar\vez$
with $\bar\vex=0$ for (\ref{Auxiliary}) with the original $\veb$ is easy
to construct. Now apply the algorithm of Lemma \ref{Iterations}
and find in time $O(n^3\bar L)=O(n^3L)$ an optimal solution $\vez$
for (\ref{Auxiliary}). If the optimal objective value is $0$ then
$\vey=\ve0$ and $\vex$ is feasible in the original program (\ref{NIP})
whereas if it is positive then (\ref{NIP}) is infeasible.
\eoproof

We can now obtain the main result of this article.

\bt{Main}
For every fixed integer $(r,s)\times t$ bimatrix $A$, there is an algorithm
that, given $n$, vectors $\vew,\vel,\veu\in\Z^{nt}$ and $\veb\in\Z^{r+ns}$
having binary encoding length $L:=\l \vew, \veb, \vel, \veu, \r$,
solves in time $O(n^3L)$ the $n$-fold integer programming problem
$$\min\left\{\vew\vex\ :\ \A\vex=\veb\,,\
\vel\leq\vex\leq\veu\,,\ \vex\in\Z^{nt}\right\}\ .$$
\et

\boproof
Use the algorithm of Lemma \ref{Feasibility} to either detect
infeasibility or obtain a feasible point and augment it by the
algorithm of Lemma \ref{Iterations} to optimality.
\eoproof

\section{Extensions to nonlinear objectives}
\label{etno}

Here we extend some of our results to programs with nonlinear objective functions,
\begin{equation}\label{SIP}
\min\left\{f(\vex)\ :\ \A\vex=\veb\,,\
\vel\leq\vex\leq\veu\,,\ \vex\in\Z^{nt}\right\}\ .
\end{equation}
A function $f:\R^{nt}\rightarrow\R$ is {\em separable convex}
if $f(\vex)=\sum_{i=1}^n f^i(\vex^i)=\sum_{i=1}^n\sum_{j=1}^t f^i_j(x^i_j)$
with each $f^i_j$ univariate convex.
In \cite{MSW} it was shown that Graver bases
provide optimality certificates for problem (\ref{SIP})
with separable convex functions as well: a feasible point $\vex$ is
optimal if and only if there is no feasible Graver step $\veg$
for $\vex$ which satisfies $f(\vex+\veg)<f(\vex)$. This was used
in \cite{HOW1} to provide polynomial time procedures for
optimality certification and solution of problem (\ref{SIP})
with separable convex functions $f$. However, this involved again
checking each of the $O(n^{g(A)})$ elements of $\G\left(\A\right)$.

Our results from Section \ref{ta} can be extended to provide
linear time optimality certification for separable convex
functions and a cubic time solution of (\ref{SIP}) for separable
convex piecewise affine functions. We discuss these respectively next.

\subsubsection*{Optimality certification for separable convex objectives}

Here we assume that the objective function $f$ is presented by a
{\em comparison oracle} that, queried on two vectors $\vex,\vey$,
asserts whether or not  $f(\vex)\leq f(\vey)$. The time complexity now
measures the number of arithmetic operations and oracle queries.

\bt{Separable}
For any fixed bimatrix $A$, there is an algorithm that, given $n$,
$\veb,\vel,\veu$, separable convex $f$ presented by comparison oracle,
and feasible point $\vex$ in program (\ref{SIP}), either asserts
that $\vex$ is optimal or finds an augmenting step $\veg$
for $\vex$ which satisfies $f(\vex+\veg)\leq f(\vex+\bar\veg)$ for
any feasible step $\bar\veg\in\G\left(\A\right)$, in linear time $O(n)$.
\et
\boproof
Given the feasible point $\vex$, set a dynamic program similar to
that in Definition \ref{DynamicProgram}, with $\gamma:=1$, with the
only modification that the weight of arc $(\veh^{i-1},\veh^i)$ from
$\veh^{i-1}\in S_{i-1}$ to $\veh^i\in S_i$ is now defined to be
$f^i(\vex^i+\veg^i)-f^i(\vex^i)$ with $\veg^i:=\veh^i-\veh^{i-1}$.
Then, for every dipath $\veh$ and its associated vector
$\veg:=\veg(\veh)$, we now have
$$f(\vex+\veg)-f(\vex)\ =\ \sum_{i=1}^n \left(f^i(\vex^i+\veg^i)-f^i(\vex^i)\right)
\ =\ \mbox{weight of dipath $\veh$}\ .$$
We now claim that the desired step is the vector $\veg:=\veg(\veh)$
associated with a minimum weight dipath $\veh$ in this dynamic program.
Indeed, an argument similar to that in the proof Lemma \ref{AugmentingStep}
now implies that for any feasible step $\bar\veg\in\G\left(\A\right)$
we have $f(\vex+\veg)-f(\vex)\leq f(\vex+\bar\veg)-f(\vex)$
and therefore $f(\vex+\veg)\leq f(\vex+\bar\veg)$.
\eoproof

\subsubsection*{Optimization of separable convex piecewise affine objectives}

In \cite{HOW1} it was shown that problem (\ref{SIP}) can be solved for
any separable convex function in polynomial time, but with very large
dependency of $O(n^{g(A)})$ of the running time on $n$, with the
exponent $g(A)$ depending on the bimatrix $A$. Here we restrict attention
to separable convex objective functions which are piecewise affine, for which
we are able to reduce the time dependency on $n$ to $O(n^3)$ independent of $A$.

So we assume again that
$f(\vex)=\sum_{i=1}^n f^i(\vex^i)=\sum_{i=1}^n\sum_{j=1}^t f^i_j(x^i_j)$
with each $f^i_j:\R\rightarrow\R$ univariate convex. Moreover,
we now also assume that for some fixed $p$, each $f^i_j$ is $p$-piecewise affine,
that is, the interval between the lower bound $\vel^i_j$ and upper bound
$\veu^i_j$ is partitioned into at most $p$ intervals with integer end-points,
and the restriction of $f^i_j$ to each interval $k$ is an affine function
$w^i_{j,k}x^i_j+a^i_{j,k}$ with all $w^i_{j,k},a^i_{j,k}$ integers.
We denote by $\l f\r$ the binary length of $f$ which is the sum
of binary lengths of all interval end-points and $w^i_{j,k},a^i_{j,k}$
needed to describe it. The binary length of the input for the nonlinear
problem (\ref{SIP}) is now $L:=\l f,\veb,\vel,\veu \r$.

\bt{Piecewise}
For any fixed $p$ and bimatrix $A$, there is an algorithm that,
given $n$, $\veb,\vel,\veu$, and separable convex $p$-piecewise
affine $f$, solves in time $O(n^3L)$ the program
$$\min\left\{f(\vex)\ :\ \A\vex=\veb\,,\
\vel\leq\vex\leq\veu\,,\ \vex\in\Z^{nt}\right\}\ .$$
\et

\boproof
We need to establish analogs of some of the lemmas
of Section \ref{ta} for such objective functions. First, for the analog
of Lemma \ref{AugmentingStep}, proceed as in the proof of
Theorem \ref{Separable} above: given a feasible point $\vex$ and positive
integer $\gamma$, set again a dynamic program similar to that in
Definition \ref{DynamicProgram}, with the weight of arc $(\veh^{i-1},\veh^i)$
from $\veh^{i-1}\in S_{i-1}$ to $\veh^i\in S_i$ defined to be
$f^i(\vex^i+\gamma\veg^i)-f^i(\vex^i)$ with $\veg^i:=\veh^i-\veh^{i-1}$.
A similar argument to that in the proof of Lemma \ref{AugmentingStep}
now shows that $\gamma\veg$ with $\veg:=\veg(\veh)$ the vector associated
with the minimum weight dipath $\veh$ is a feasible step for $\vex$ which
satisfies $f(\vex+\gamma\veg)\leq f(\vex+\gamma\bar\veg)$ for any
feasible step $\gamma\bar\veg$ with $\bar\veg\in\G\left(\A\right)$.

For the analog of Lemma \ref{BestGraverStep}, we construct again a
set $\Gamma$ of critical step sizes $\gamma$ as follows. First, as in
the proof of Lemma \ref{BestGraverStep}, we collect the critical step
sizes due to the lower and upper bound constraints by finding,
for every $i=1,\dots,n$ and every $\vez\in Z(A)\setminus\{\ve 0\}$,
the largest positive integer $\gamma$ such that
$\vel^i\leq\vex^i+\gamma\vez\leq\veu^i$, and include it in $\Gamma$.
However, in contrast to Lemma \ref{BestGraverStep}, we are now dealing with
the more general class of piecewise affine objective functions $f^i_j$.
So we must add also the following values $\gamma$ to $\Gamma$: for every
$i=1,\dots,n$, every $\vez\in Z(A)\setminus\{\ve 0\}$, and every $j=1,\ldots,t$,
if $x^i_j+\gamma z_j$ and $x^i_j+(\gamma+1) z_j$ belong to different affine
pieces of $f^i_j$, then $\gamma$ is included in $\Gamma$. Since the number
$p$ of affine pieces is constant, the number of such values for each $i$,
$z$ and $j$ is also constant. So the total number of elements of $\Gamma$
remains linear and it can be constructed in linear time $O(n)$ again.
Now we continue as in the proof of Lemma \ref{BestGraverStep}:
for each $\gamma$ in $\Gamma$, using the analog of Lemma \ref{AugmentingStep}
established in the first paragraph above, we solve the corresponding
dynamic program in $O(n)$ time, resulting in total of $O(n^2)$ time again.
Let $\gamma\veg$ be that feasible step for $\vex$ which attains minimum
value $f(\vex+\gamma\veg)$ among the best steps obtained from
all these dynamic programs, and let $\bar\gamma\bar\veg$ be that
feasible step for $\vex$ which attains minimum value
$f(\vex+\bar\gamma\bar\veg)$ among $\bar\gamma\bar\veg$
with $\bar\veg\in\G\left(\A\right)$ if any. It now follows from the
construction of $\Gamma$ that if $\bar\gamma\bar\veg$ is augmenting,
namely, if $f(\vex+\bar\gamma\bar\veg)-f(\vex)<0$, then $\bar\gamma\in\Gamma$,
as otherwise the step $(\bar\gamma+1)\bar\veg$ will be
feasible and better. Let $\bar\gamma\hat\veg$ be the best
step attained from the dynamic program of $\bar\gamma$. Then
$$f(\vex+\gamma\veg)\ \leq\ f(\vex+\bar\gamma\hat\veg)
\ \leq\ f(\vex+\bar\gamma\bar\veg)$$
where the first inequality follows from the choice of $\gamma\veg$
and the second inequality follows from the analog of
Lemma \ref{AugmentingStep}. Therefore
$f(\vex+\gamma\veg)\leq f(\vex+\bar\gamma\bar\veg)$.

For the analog of Lemma \ref{Iterations}, we use the results of
\cite{HOW1} incorporating the optimality criterion of \cite{MSW},
which assure that the number of iterations needed when using a
Graver-best augmenting step at each iteration, is bounded by
$O(n\,\l f \r)=O(nL)$, resulting again in overall time complexity
$O(n^3 L)$ for augmenting an initial feasible point to an optimal
solution of (\ref{SIP}). Since an initial feasible point if any
can be found by Lemma \ref{Feasibility} as before in the same
complexity, the theorem now follows.
\eoproof

\section{Some consequences}
\label{sc}

Here we briefly discuss two of the many consequences
of $n$-fold integer programming which, now with our
new algorithm, can be solved drastically faster than before.

\subsubsection*{Nonlinear multicommodity transportation}

The multicommodity transportation problem seeks minimum cost
routing of $l$ commodities from $m$ suppliers to $n$
consumers subject to supply, consumption and capacity constraints.
For $l=1$ this is the classical transportation
problem which is efficiently solvable by linear programming.
But already for $l=2$ it is NP-hard. Here we consider the problem
with fixed (but arbitrary) number $l$ of commodities, fixed (but arbitrary) number
$m$ of suppliers, and variable number $n$ of consumers. This is natural
in typical applications where few facilities serve many customers.

The data is as follows. Each supplier $i$ has a supply vector $s^i\in\Z_+^l$
with $s^i_k$ its supply in commodity $k$. Each consumer $j$ has a consumption
vector $c^j\in\Z_+^l$ with $c^j_k$ its consumption in commodity $k$.
The amount of commodity $k$ to be routed from supplier $i$ to consumer $j$
is an integer decision variable $x^j_{i,k}$.
The total amount $\sum_{k=1}^l x^j_{i,k}$ of commodities routed on
the channel from $i$ to $j$ should not exceed the channel capacity $u_{i,j}$,
and has cost $f_{i,j}\left(\sum_{k=1}^l x^j_{i,k}\right)$ for suitable
univariate functions $f_{i,j}$. We can handle standard linear costs as
well as more realistic, convex piecewise affine cost functions $f_{i,j}$,
which account for channel congestion under heavy routing.

As a corollary of Theorem \ref{Piecewise}, for any
fixed numbers $l$ of commodities and $m$ of suppliers,
the problem can be solved in time cubic in the number $n$ of consumers.

\bc{Multicommodity}
For every fixed $l$ commodities, $m$ suppliers, and $p$,
there exists an algorithm that, given $n$ consumers,
supplies and demands $s^i,c^j\in\Z_+^l$, capacities $u_{i,j}\in\Z_+$,
and convex $p$-piecewise affine costs $f_{i,j}:\Z\rightarrow\Z$,
solves in time $O(n^3 L)$, with $L:=\l s^i,c^j,u_{i,j},f_{i,j}\r$,
the integer multicommodity transportation problem
\begin{eqnarray*}
\min\left\{\sum_{i=1}^m\sum_{j=1}^n f_{i,j}\left(\sum_{k=1}^l x^j_{i,k}\right)\, :\,
x^j_{i,k}\in\Z_+,\ \sum_j x^j_{i,k}=s^i_k,\ \sum_i x^j_{i,k}=c^j_k,\
\sum_{k=1}^l x^j_{i,k}\leq u_{i,j}\right\}\,.
\end{eqnarray*}
\ec

\boproof
Introduce new variables $y^j_i$ and equations
$y^j_i=\sum_{k=1}^l x^j_{i,k}$ for all $i,j$.
Then the objective function becomes $\sum_{i,j}f_{i,j}(y^j_i)$
which is separable convex $p$-piecewise affine in the new variables,
and the capacity constraints become $y^j_i\leq u_{i,j}$ which
are upper bounds on the new variables. Use $u_{i,j}$ as an upper
bound on $x^j_{i,k}$ and $0$ as a trivial
lower bound on $y^j_i$ for all $i,j,k$.
As shown in \cite{HOW2},
arranging the original and new variables in a tuple
$\vez=(\vez^1,\dots,\vez^n)$ of $n$ bricks
$\vez^j\in\Z^{m(l+1)}$ defined by
$$\vez^j:\ =\ (x^j_{1,1},\dots,x^j_{1,l},y^j_1\,,\
x^j_{2,1},\dots,x^j_{2,l},y^j_2\,,\ \dots\dots\,,\
x^j_{m,1},\dots,x^j_{m,l},y^j_m)\ ,$$
this problem can be modeled as an $n$-fold program, resulting in solution
in large time $O\left(n^{g(A)}L\right)$, with exponent depending on $l,m$.
Theorem \ref{Piecewise} now enables solution in cubic time
independent of the numbers $l$ of commodities and $m$ of suppliers.
\eoproof

\subsubsection*{Privacy in statistical databases}

A common practice in the disclosure of sensitive data contained in a
multiway table is to release some of the table margins rather than the
entries of the table. Once the margins are released, the security of any
specific entry of the table is related to the set of possible values that can
occur in that entry in all tables having the same margins as those of the
source table in the database, see \cite{DFRSZ,SZP} and the references
therein. In particular, if this set is small or consists of a unique value,
that of the source table, then this entry can be exposed. Thus, it
is desirable to compute the minimum and maximum integer values that
can occur in an entry, which in particular are equal if and only
if the entry value is unique, before margin disclosure is enabled.

Consider $(d+1)$-way tables of format $m_0\times\cdots\times m_d$, that
is, arrays $v=(v_{i_0,\dots,i_d})$ indexed by $1\leq i_j\leq m_j$
for all $j$, with all entries $v_{i_0,\dots,i_d}$ nonnegative integers.
Our results hold for arbitrary hierarchical margins, but for simplicity
we restrict attention to disclosure of $d$-margins, that is, the $d+1$ many
$d$-way tables $(v_{i_0,\dots,i_{j-1},*,i_{j+1},\dots, i_d})$ obtained from
$v$ by collapsing one factor $0\leq j\leq d$ at a time, with entries given by
$$v_{i_0,\dots,i_{j-1},*,i_{j+1},\dots, i_d}:=
\sum_{i_j=1}^{m_j}v_{i_0,\dots,i_{j-1},i_j,i_{j+1},\dots, i_d}\,,\quad
1\leq i_k\leq m_k\,,\ \ 0\leq k\leq d\,,\ \ k\neq j\ \ .$$
The problem is then to compute the minimum and maximum integer
values that can occur in an entry subject to the margins
of the source table in the database.

This problem is NP-hard
already for $3$-way tables of format $n\times m\times 3$, see \cite{DO}.
However, as a corollary of Theorem \ref{Main}, if only one
side $n$ of the table is variable, the problem can be solved
in cubic time regardless of the other sides $m_i$ as follows.

\bc{Privacy}
For every fixed $d,m_1,\dots,m_d$, there is an algorithm that, given $n$,
integer $d$-margins $(v_{*,i_1,\dots, i_d}),\dots,(v_{i_0,\dots,i_{d-1},*})$,
and index $(k_0,\dots,k_d)$, determines, in time $O(n^3 L)$,
with $L$ the binary length of the given margins, the minimum and
maximum values of entry $x_{k_0,\dots,k_d}$ among all tables
with these margins, that is, solves
$$\min/\max\left\{x_{k_0,\dots,k_d}:
x\in\Z_+^{n\times m_1\times\cdots\times m_d},\,
(x_{i_0,\dots,i_{j-1},*,i_{j+1},\dots, i_d})=
(v_{i_0,\dots,i_{j-1},*,i_{j+1},\dots, i_d})\ \forall j\right\}.$$
\ec

\boproof
Let $u$ be the maximum value of any entry in the given margins
and use it as an upper bound on every variable.
As shown in \cite{DHOW}, this problem can be modeled as an
$n$-fold integer programming problem, resulting in solution
in large running time $O\left(n^{g(A)}L\right)$,
with exponent which depends on $m_1,\dots,m_d$.
Theorem \ref{Main} now enables to solve it in cubic time
independent of the table dimensions $m_1,\dots,m_d$.
\eoproof

We note that long tables, with one side much larger
than the others, often arise in practical applications.
For instance, in health statistical tables, the long factor may be
the age of an individual, whereas other factors may be binary
(yes-no) or ternary (subnormal, normal, and supnormal).
Moreover, it is always possible to merge categories of factors,
with the resulting coarser tables approximating the original ones,
making the algorithm of Corollary \ref{Privacy} applicable.

We also note that, by repeatedly incrementing a lower bound
and decrementing an upper bound on the entry $x_{k_0,\dots,k_d}$,
and computing its new minimum and maximum values subject to these
bounds, we can produce the entire set of values that can occur in
that entry in time proportional to the number of such values.

\section{Solvability over bimatrices with variable entries}
\label{sobwve}

The drop of the Graver complexity from the exponent of $n$ to
the constant multiple also leads to the first polynomial
time solution of $n$-fold integer programming with variable
bimatrices. Of course, by the universality of $n$-fold integer
programming, the variability of the bimatrices must be limited.
In what follows, we fix the dimensions $r,s,t$ of the input
bimatrix $A$, and let the entries vary. We show that, given
as part of the input an upper bound $a$ on the absolute value
of every entry of $A$, we can solve the problem in time polynomial
in $a$, that is, polynomial in the unary length of $a$. This holds
for linear as well as separable convex piecewise affine objectives.

We have the following theorem,
with $L:=\l f,a,\veb,\vel,\veu \r$ the length of the input.

\bt{Variable}
For any fixed $r,s,t$ and $p$, there is an algorithm that,
given $n$, $a$, $(r,s)\times t$ bimatrix $A$ with all entries
bounded by $a$ in absolute value, $\veb,\vel,\veu$, and separable
convex $p$-piecewise affine $f$, in polynomial
time $O(a^{3t(rs+st+r+s)}n^3L)$, solves
$$\min\left\{f(\vex)\ :\ \A\vex=\veb\,,\
\vel\leq\vex\leq\veu\,,\ \vex\in\Z^{nt}\right\}\ .$$
\et

\boproof
Let $\G(A_2)$ be the Graver basis of the $s\times t$ second block
$A_2$ of $A$, let $p:=|\G(A_2)|$ be its cardinality, and arrange
its elements as the columns of a $t\times p$ matrix $G_2$.
Since $r,s,t$ are fixed, it follows from bounds on Graver bases
(see e.g. \cite[Section 3.4]{Onn}) that every $g\in\G(A_2)$
satisfies $\|g\|_\infty=O(a^s)$ and hence $p=O(a^{st})$.

Now, it is known (see \cite{HS,SS} or \cite[Section 4.1]{Onn})
that the Graver complexity $g(A)$ of $A$ is equal to the
maximum value $\|v\|_1$ of any element $v$ in the Graver
basis $\G(A_1G_2)$ of the $r\times p$ matrix $A_1G_2$.
Since the entries of $A_1G_2$ are bounded in absolute value by
$O(a^{s+1})$, the bounds on Graver bases
(see again \cite[Section 3.4]{Onn}) imply that
$\|v\|_1=O(p\cdot(a^{s+1})^r)$ for every $v\in\G(A_1G_2)$
and hence $g(A)=O(a^{rs+st+r})$.

Now, consider again the following set
defined in (\ref{Space}) in Lemma \ref{StateSpace},
$$Z(A)\ :=\ \left\{\vez\in\Z^t\,:\, \vez\
\mbox{is the sum of at most $g(A)$ elements of $\G(A_2)$}\right\}\ .$$
For each $z\in Z(A)$ we have that $\|z\|_\infty$ is bounded
by $g(A)$ times the maximum value of $\|g\|_\infty$ over all
$g\in\G(A_2)$, and therefore
$\|z\|_\infty=O(g(A)\cdot a^s)=O(a^{rs+st+r+s})$.
So the cardinality of $Z(A)$ satisfies
$|Z(A)|=O((a^{rs+st+r+s})^t)=O(a^{t(rs+st+r+s)})$.

Now, suitable analogs of some of the lemmas of Section \ref{ta}
go through, except that the complexities now depend on the
variable size of the set $Z(A)$, as follows. The time complexity of
the algorithm of Lemma \ref{AugmentingStep} becomes $O(|Z(A)|^2n)$.
The size of the set $\Gamma$ of critical values is now
$O(|Z(A)|n)$ and therefore the time complexity of the
algorithm of Lemma \ref{BestGraverStep} now becomes $O(|Z(A)|^3n^2)$.
The number of iterations needed to augment an initial feasible
point to an optimal solution remains $O(nL)$ as before and therefore
the time complexity of the algorithm of Lemma \ref{Iterations}
now become $O(|Z(A)|^3n^3L)$. To find an initial feasible point,
one can use the algorithm of Lemma \ref{Feasibility}, but then
$t$ would have to be replaced by $t+2r+2s$ for the auxiliary
bimatrix $\bar A$, resulting in a somewhat larger exponent for $a$
in the running time. However, it is possible to find an
initial feasible point in an alternative, somewhat more involved way,
keeping the original system with the bimatrix $A$, as follows.
First find an integer solution to the system of equations only
(without the lower and upper bounds) using the Hermite normal
forms of the blocks $A_1$ and $A_2$. Second, relax the bounds so
as to make that point feasible. Third, minimize the following auxiliary
objective function which is separable convex $3$-piecewise affine, with
$$f^i_j(x^i_j)\ :=\ \left\{\begin{array}{ll}
      l^i_j-x^i_j, & \text{if } x^i_j\leq l^i_j,\\
      0, & \text{if } l^i_j\leq x^i_j\leq u^i_j,\\
      x^i_j-u^i_j, & \text{if } x^i_j\geq u^i_j.
    \end{array}\right.$$
If the optimal value is zero then the optimal auxiliary solution
is feasible in the original problem, whereas if it is positive
then the original problem is infeasible. Since this minimization
can be done using the separable convex piecewise affine analog of
Lemma \ref{Iterations} described in the proof of
Theorem \ref{Piecewise} in the same complexity $O(|Z(A)|^3n^3L)$,
the overall running time is $O(a^{3t(rs+st+r+s)}n^3L)$ as claimed.
\eoproof

\section{Parametrization and approximation hierarchy}
\label{paah}

We conclude with a short discussion of the universality of $n$-fold
integer programming and the resulting parametrization and simple
approximation hierarchy for all of integer programming. As mentioned
in the introduction, every integer program is an $n$-fold program
for some $m$ over the bimatrix $A(m)$ having first block the identity
matrix $I_{3m}$
and second block the $(3+m)\times 3m$ incidence matrix of $K_{3,m}$. It is
convenient and illuminating to introduce also the following description.

Consider the following special form of the $n$-fold product operator.
For an $s\times t$ matrix $D$, let $D^{[n]}:=A^{(n)}$
where $A$ is the $(t,s)\times t$ bimatrix $A$ with first block
$A_1:=I_t$ the $t\times t$ identity matrix and second block $A_2:=D$.
We consider such $m$-fold products of the $1\times 3$ matrix $(1\ 1\ 1)$.
Note that $(1\ 1\ 1)^{[m]}$ is precisely the $(3+m)\times 3m$
incidence matrix of the complete bipartite graph $K_{3,m}$. For instance,
\begin{equation*}\label{matrix}
(1\ 1\ 1)^{[2]}\ =\
\left(
\begin{array}{cccccc}
  1 & 0 & 0 & 1 & 0 & 0 \\
  0 & 1 & 0 & 0 & 1 & 0 \\
  0 & 0 & 1 & 0 & 0 & 1 \\
  1 & 1 & 1 & 0 & 0 & 0 \\
  0 & 0 & 0 & 1 & 1 & 1 \\
\end{array}
\right)\ .
\end{equation*}

The following theorem was established in \cite{DO}.

\vskip.2cm\noindent
{\bf The Universality Theorem \cite{DO}}
{\em Every (bounded) integer programming problem
$\min\{\vecc\vey\,:\,\vey\in\Z_+^k,\ V\vey=\vev\}$
is polynomial time equivalent to some integer program
\begin{equation*}\label{Universal}
\min\left\{\vew\vex\, :\, \vex\in\Z_+^{3mn},\, (1\ 1\ 1)^{[m][n]}\vex=\veb\right\}
\ \cong\ \min\left\{\vew\vex\, :\, \vex\in\Z_+^{3mn},\, A(m)^{(n)}\vex=\veb\right\}\,.
\end{equation*}}

This theorem provides a new, variable dimension, parametrization
of integer programming: for each fixed value of the parameter $m$,
the resulting programs above with $n$ variable live in
variable dimension $3mn$ and include natural models such as
those described in Section \ref{sc}, and can be solved in cubic time $O(n^3L)$
by Theorem \ref{Main}; and when the parameter $m$ varies,
every integer program appears for some $m$.

Our new algorithm suggests a natural simple approximation hierarchy
for integer programming, parameterized by degree $d$, as follows.
Fix any $d$. Then given any $m$, let $A:=A(m)$, so $t=3m$,
$A_1=I_{3m}$, and $A_2$ is the incidence matrix of $K_{3,m}$.
Define the approximation at degree $d$ of the set $Z(A)$ in
equation (\ref{Space}) in Lemma \ref{StateSpace} by
\begin{eqnarray}\label{ApproximatedSpace}
Z_d(m)\ :=\ \left\{\vez\in\Z^{3m}\,:\, \vez\
\mbox{is the sum of at most $d$ elements of $\G(A_2)$}\right\}\ .
\end{eqnarray}
Since $A_2$ is totally unimodular, $\G(A_2)$ consists of the
$O(m^3)$ vectors in $\{0,\pm1\}^{3m}$ supported on circuits of
$K_{3,m}$ with alternating $\pm1$ and hence $|Z_d(m)|=O(m^{3d})$.

Now, given a feasible point $\vex$ in the universal program
above, and positive integer $\gamma$, set a dynamic program similar
to that in Definition \ref{DynamicProgram}, with the only modification that
the sets $S_i$ are defined using the approximation $Z_d(m)$ of $Z(A)$.
Weaker forms of Lemmas \ref{AugmentingStep} and \ref{BestGraverStep}
now assert that in time $O(|Z_d(m)|^3 n^2)=O(m^{9d}n^2)$, which is
polynomial in both $m$ and $n$, we can find a good feasible step $\gamma\veg$.
We use this iteratively to augment an initial feasible point to
one which is as good as possible and output it. However, Lemma \ref{StateSpace}
no longer holds and not all brick sums of elements of the Graver basis $\G(\A)$
lie in $Z_d(m)$. So the bounds on the number of iterations and total
running time are no longer valid and the output point may be non optimal.

By increasing the degree $d$ we can get better approximations at
increasing running times, and when $d=g(A)$ we get the true optimal
solution. These approximations are currently under study, implementation
and testing. They show promising behavior already at degree $d=3$
and will be discussed in more detail elsewhere.

For $m=3$, discussed in Example \ref{Example}, for which the universal problem
is equivalent to optimization over $3$-way $n\times 3\times 3$ tables,
the approximation $Z_3(3)$ at degree $d=3$ contains only $811$ vectors
out of the $42931$ vectors in the true $Z(A)$, such as
$$\left(\begin{array}{ccccccccc}
-3 & 2 & 1 & 2 & -3 & 1 & 1 & 1 & -2 \end{array}\right)\ .$$


\begin{thebibliography}{}

\bibitem{AT}
Aoki, S., Takemura, A.:
Minimal basis for connected Markov chain over $3\times3\times K$
contingency tables with fixed two-dimensional marginals.
Austr. New Zeal. J. Stat. 45 (2003) 229--249

\bibitem{BO}
Berstein, Y., Onn, S.:
The Graver complexity of integer programming.
Ann. Combin. 13 (2009) 289--296

\bibitem{CFS}
Cook, W., Fonlupt, J., Schrijver, A.:
An integer analogue of Carath\'eodory's theorem.
J. Combin. Theory Ser. B 40 (1986) 63--70

\bibitem{DHORW}
De Loera, J., Hemmecke, R., Onn, S., Rothblum, U.G., Weismantel, R.:
Convex integer maximization via Graver bases.
J. Pure Appl. Algebra. 213 (2009) 1569--1577

\bibitem{DHOW}
De Loera, J.A., Hemmecke, R., Onn, S., Weismantel, R.:
N-fold integer programming.
Discrete Optimization 5 (Volume in memory of George B. Dantzig),
231--241 (2008)

\bibitem{DO}
De Loera, J.A., Onn, S.:
All linear and integer programs are slim 3-way transportation programs.
SIAM Journal on Optimization 17, 806--821 (2006)

\bibitem{DFRSZ}
Dobra, A., Fienberg, S.E., Rinaldo, A., Slavkovi\'c, A., Zhou, Y.:
Algebraic statistics and contingency table problems: log-linear
models, likelihood estimation, and disclosure limitation.
In: Emerging Applications of Algebraic Geometry: IMA Volumes
in Mathematics and its Applications 148 (2009) 63--88, Springer

\bibitem{Gra}
Graver, J.E.:
On the foundation of linear and integer programming {I}.
Mathematical Programming 9, 207--226 (1975)

\bibitem{HKW}
Hemmecke, R., K\"oppe, M., Weismantel, R.:
A polynomial-time algorithm for optimizing over
N-fold 4-block decomposable integer programs. IPCO 14 (2010)

\bibitem{HOW1}
Hemmecke, R., Onn, S., Weismantel, R.:
A polynomial oracle-time algorithm for convex integer minimization.
Mathematical Programming 126, 97--117, (2011)

\bibitem{HOW2}
Hemmecke, R., Onn, S., Weismantel, R.:
N-fold integer programming and nonlinear multi-transshipment.
Optimization Letters, 5, 13--25 (2011)

\bibitem{HSc}
Hemmecke, R., Schultz, R.:
Decomposition of test sets in stochastic integer programming.
Mathematical Programming, 94, 323–-341 (2003)

\bibitem{HS}
Ho\c sten, S., Sullivant, S.:
Finiteness theorems for Markov bases of hierarchical models.
Journal of Combinatorial Theory, Series A, 114, 311--321 (2007)

\bibitem{KMW}
Kobayashi, Y., Murota, K., Weismantel, R.:
Cone superadditivity of discrete convex functions.
METR 2009-30, University of Tokyo, Japan, 25 pp. (2009)

\bibitem{LS}
Louveaux, F.V., Schultz, R.:
Stochastic Integer Programming.
In: Handbooks in Operations Research
and Management Science 10 (2003) 213–-266, Elsevier

\bibitem{MSW}
Murota, K., Saito, H., Weismantel, R.:
Optimality criterion for a class of nonlinear integer programs.
Operations Research Letters 32, 468--472 (2004)

\bibitem{Onn}
Onn, S.: {\bf Nonlinear Discrete Optimization}.
Zurich Lectures in Advanced Mathematics,
European Mathematical Society, x+137 pp. (September 2010)

\bibitem{SS}
Santos, F., Sturmfels, B.:
Higher Lawrence configurations.
Journal of Combinatorial Theory, Series A, 103, 151--164 (2003)

\bibitem{Seb}
Seb\"o, A.:
Hilbert bases, Carath\'eodory's theorem
and combinatorial optimization.
Waterloo University Press, IPCO 1, 431--455 (1990)

\bibitem{SZP}
Slavkovi\'c, A.B., Zhu, X., Petrovi\'c, S.:
A sample space of k-way tables given conditionals and their relations
to marginals: Implications for cell bounds and Markov bases.
Preprint, 35 pp. (2009)

\end{thebibliography}
\end{document}